\documentclass[12pt,a4paper,twoside]{article}

\newcommand{\pageformat}[6]{\setlength{\hoffset}{-1in}
                  \setlength{\voffset}{-1in}
                  \addtolength{\hoffset}{#5}
                            \addtolength{\voffset}{#6}
                            \setlength{\oddsidemargin}{#1}
                            \setlength{\evensidemargin}{#2}
                            \setlength{\textwidth}{\paperwidth}
                  \addtolength{\textwidth}{-\oddsidemargin}
                  \addtolength{\textwidth}{-\evensidemargin}
                  \addtolength{\textwidth}{-\marginparsep}
                  \addtolength{\textwidth}{-\marginparwidth}
                            \setlength{\topmargin}{#3}
                            \setlength{\textheight}{\paperheight}
                  \addtolength{\textheight}{-\topmargin}
                  \addtolength{\textheight}{-\headheight}
                  \addtolength{\textheight}{-\headsep}
                  \addtolength{\textheight}{-\footskip}
                  \addtolength{\textheight}{-#4}}
\pageformat{2cm}{3cm}{25mm}{25mm}{1pt}{0pt}

\usepackage{ifthen}
\newboolean{article}
    \setboolean{article}{true}
\newboolean{report}
\newboolean{book}
\newboolean{letter}
\newboolean{german}
\newboolean{italian}
\newboolean{nobaselinestretch}
\newboolean{nosectionappendix}
\newboolean{oldtoc}
\newboolean{nosectionequn}
\newboolean{notheorem}
\ifthenelse{\boolean{german}}{
    \usepackage{german}}{}

\usepackage[latin1]{inputenc}
\usepackage{amsmath}
\usepackage{amssymb}
\usepackage[mathscr]{eucal}

\ifthenelse{\boolean{notheorem}}{}{
    \usepackage{theorem}}

\ifthenelse{\boolean{nobaselinestretch}}{}{
    \renewcommand{\baselinestretch}{1.25}}

\newenvironment{env}[2]{\begin{#1}#2\end{#1}}{}
    \newcommand{\beq}[1]{\begin{env}{equation}{#1}}
    \newcommand{\beqn}[1]{\begin{env}{equation*}{#1}}
    \newcommand{\bal}[1]{\begin{env}{align}{#1}}
    \newcommand{\baln}[1]{\begin{env}{align*}{#1}}
    \newcommand{\bga}[1]{\begin{env}{gather}{#1}}
    \newcommand{\bgan}[1]{\begin{env}{gather*}{#1}}
    \newcommand{\bflal}[1]{\begin{env}{flalign}{#1}}
    \newcommand{\bflaln}[1]{\begin{env}{flalign*}{#1}}
    \newcommand{\bmu}[1]{\begin{env}{multline}{#1}}
    \newcommand{\bmun}[1]{\begin{env}{multline*}{#1}}
    \newcommand{\bsp}[1]{\begin{env}{split}{#1}}

    \newcommand{\eeq}{\end{env}}
    \newcommand{\eeqn}{\end{env}}
    \newcommand{\eal}{\end{env}}
    \newcommand{\ealn}{\end{env}}
    \newcommand{\ega}{\end{env}}
    \newcommand{\egan}{\end{env}}
    \newcommand{\eflal}{\end{env}}
    \newcommand{\eflaln}{\end{env}}
    \newcommand{\emu}{\end{env}}
    \newcommand{\emun}{\end{env}}
    \newcommand{\esp}{\end{env}}

\newcommand{\lf}{\vspace{2ex}}

\renewcommand{\bf}[1]{\textbf{#1}}
\renewcommand{\it}[1]{\textit{#1}}

\renewcommand{\sf}[1]{\textsf{#1}}

\renewcommand{\tt}[1]{\texttt{#1}}
\newcommand{\hl}[1]{\bf{\it{#1}}}

\newcommand{\mbf}[1]{\mathbf{#1}}
\newcommand{\msf}[1]{\text{\small$\sf{#1}$}}

\newcommand{\cmc}[1]{\mathcal{#1}}
\newcommand{\eus}[1]{\mathscr{#1}}

\newcommand{\bb}[1]{\mathbb{#1}}

\newcommand{\nbd}[1]{$#1$\nobreakdash--}
\newcommand{\ol}[1]{\overline{#1}}

\newcommand{\wt}[1]{\widetilde{#1}}

\newcommand{\vp}{\varphi}

\newcommand{\bfam}[1]{\bigl(#1\bigr)}

\newcommand{\AB}[1]{\langle#1\rangle}

\newcommand{\CB}[1]{\{#1\}}
\newcommand{\bCB}[1]{\bigl\{#1\bigr\}}
\newcommand{\BCB}[1]{\Bigl\{#1\Bigr\}}

\newcommand{\set}[2][]{
    \ifthenelse{\equal{#1}{}}{
        \CB{#2}}{
        \CB{#1~|~#2}}}
\newcommand{\bset}[2][]{
    \ifthenelse{\equal{#1}{}}{
        \bCB{#2}}{
        \bCB{#1~|~#2}}}
\newcommand{\Bset}[2][]{
    \ifthenelse{\equal{#1}{}}{
        \BCB{#2}}{
        \BCB{#1~\big|~#2}}}

\DeclareMathOperator{\ls}{\normalfont\msf{span}}
\DeclareMathOperator{\cls}{\ol{\ls}}

\DeclareMathOperator{\id}{\normalfont\msf{id}}

\newcommand{\C}{\bb{C}}

\newcommand{\cA}{\cmc{A}}
\newcommand{\cB}{\cmc{B}}
\newcommand{\cC}{\cmc{C}}

\newcommand{\sB}{\eus{B}}

\newcommand{\sE}{\eus{E}}
\newcommand{\sF}{\eus{F}}

\newcommand{\sK}{\eus{K}}

\newcommand{\U}{\mbf{1}}

\ifthenelse{\boolean{nosectionequn}}{}{
    \numberwithin{equation}{section}
    }

\ifthenelse{\boolean{article}\or\boolean{letter}\or\boolean{nosectionequn}}{
    \setboolean{nosectionappendix}{true}}{}
\ifthenelse{\boolean{nosectionappendix}}{}{
    \renewcommand{\appendix}{
        \chapter*{\appendixname}
        \addcontentsline{toc}{chapter}{\appendixname}
        \renewcommand{\thesection}{\Alph{section}}
        \setcounter{section}{0}}}
   
\ifthenelse{\boolean{report}\or\boolean{book}}{
    }{}

\ifthenelse{\boolean{notheorem}}{}{
        \newcommand{\mnname}{Mathematical note.}
        \newcommand{\enname}{End of the note.}
        \newcommand{\definame}{Definition.}
        \newcommand{\propname}{Proposition.}
        \newcommand{\lemname}{Lemma.}
        \newcommand{\exname}{Example.}
        \newcommand{\exername}{Exercise.}
        \newcommand{\remname}{Remark.}
        \newcommand{\obname}{Observation.}
        \newcommand{\thmname}{Theorem.}
        \newcommand{\corname}{Corollary.}
        \newcommand{\proofname}{Proof.}
    \ifthenelse{\boolean{german}}{
        \renewcommand{\mnname}{Mathematische Notiz.}
        \renewcommand{\enname}{Ende der Notiz.}
        \renewcommand{\exname}{Beispiel.}
        \renewcommand{\exername}{Übung.}
        \renewcommand{\remname}{Bemerkung.}
        \renewcommand{\obname}{Beobachtung.}
        \renewcommand{\thmname}{Satz.}
        \renewcommand{\corname}{Korollar.}
        \renewcommand{\proofname}{Beweis.}}{}
    \ifthenelse{\boolean{italian}}{
        \renewcommand{\mnname}{Nota matematica.}
        \renewcommand{\enname}{Fina della nota.}
        \renewcommand{\definame}{Definizione.}
        \renewcommand{\propname}{Proposizione.}
        \renewcommand{\exname}{Esempio.}
        \renewcommand{\exername}{Esercizio.}
        \renewcommand{\remname}{Nota.}
        \renewcommand{\obname}{Osservazione.}
        \renewcommand{\thmname}{Teorema.}
        \renewcommand{\corname}{Corollario.}
        \renewcommand{\proofname}{Dimostrazione.}

       \renewcommand{\appendixname}{Appendice}

       }{}
    \theoremheaderfont{\normalfont\bfseries}
    \theoremstyle{change}
        \theorembodyfont{\rmfamily}
            \newtheorem{emp}{}[section]
                \newcommand{\bemp}[1][]{
                    \begin{emp}\hskip-\labelsep\bf{#1}\hskip\labelsep}
                \newcommand{\eemp}{\end{emp}}
\newtheorem{itemp}[emp]{}
                \newcommand{\bitemp}[1][]{
                    \begin{itemp}\hskip-\labelsep\bf{#1}\hskip\labelsep\normalfont\itshape}
                \newcommand{\eitemp}{\end{itemp}}
            \newtheorem{mn}[emp]{\mnname}
                \newcommand{\bnm}{\begin{mn}~\begin{quotation}\renewcommand{\baselinestretch}{1}\small\noindent\ignorespaces}
                \newcommand{\enm}{\end{quotation}\hfill\bf{\enname}\end{mn}}
            \newtheorem{ex}[emp]{\exname}
                \newcommand{\bex}{\begin{ex}}
                \newcommand{\eex}{\end{ex}}
            \newtheorem{exer}[emp]{\exername}
                \newcommand{\bexer}{\begin{exer}}
                \newcommand{\eexer}{\end{exer}}
            \newtheorem{defi}[emp]{\definame}
                \newcommand{\bdefi}{\begin{defi}}
                \newcommand{\edefi}{\end{defi}}
            \newtheorem{rem}[emp]{\remname}
                \newcommand{\brem}{\begin{rem}}
                \newcommand{\erem}{\end{rem}}
            \newtheorem{ob}[emp]{\obname}
                \newcommand{\bob}{\begin{ob}}
                \newcommand{\eob}{\end{ob}}
        \theorembodyfont{\normalfont\itshape}
            \newtheorem{thm}[emp]{\thmname}
                \newcommand{\bthm}{\begin{thm}}
                \newcommand{\ethm}{\end{thm}}
            \newtheorem{prop}[emp]{\propname}
                \newcommand{\bprop}{\begin{prop}}
                \newcommand{\eprop}{\end{prop}}
            \newtheorem{cor}[emp]{\corname}
                \newcommand{\bcor}{\begin{cor}}
                \newcommand{\ecor}{\end{cor}}
            \newtheorem{lem}[emp]{\lemname}
                \newcommand{\blem}{\begin{lem}}
                \newcommand{\elem}{\end{lem}}
\newenvironment{empn}[1]{\lf\noindent\bf{#1}\ignorespaces\hskip\labelsep}{\lf}
		\newcommand{\bempn}[1]{\begin{empn}{#1}}
		\newcommand{\eempn}{\end{empn}}
		\newcommand{\bitempn}[1]{\begin{empn}{#1}\normalfont\itshape}
		\newcommand{\eitempn}{\end{empn}}
                \newcommand{\bnmn}{\begin{empn}{\mnname}~\begin{quotation}\renewcommand{\baselinestretch}{1}\small\noindent\ignorespaces}
                \newcommand{\enmn}{\end{quotation}\hfill\bf{\enname}\end{empn}}
		\newcommand{\bexn}{\begin{empn}{\exname}}
		\newcommand{\eexn}{\end{empn}}
		\newcommand{\bexern}{\begin{empn}{\exername}}
		\newcommand{\eexern}{\end{empn}}
		\newcommand{\bdefin}{\begin{empn}{\definame}}
		\newcommand{\edefin}{\end{empn}}
		\newcommand{\bremn}{\begin{empn}{\remname}}
		\newcommand{\eremn}{\end{empn}}
		\newcommand{\bobn}{\begin{empn}{\obname}}
		\newcommand{\eobn}{\end{empn}}
		\newcommand{\bthmn}{\bitempn{\thmname}}
		\newcommand{\ethmn}{\eitempn}

		\newcommand{\bcorn}{\bitempn{\corname}}
		\newcommand{\ecorn}{\eitempn}

\newcommand{\qedsymbol}{~\rule[-0.35mm]{2mm}{2mm}}
    \newcounter{proof}[emp]
    \newenvironment{Proof}[1]{
        \vspace{1ex}
        \renewcommand{\item}[1][\stepcounter{proof}(\roman{proof})]%
            {##1\hskip\labelsep}
        \noindent\textsc{#1\hskip\labelsep}}{
        \nolinebreak\qedsymbol}
    \newcommand{\proof}[1][\proofname]{
        \begin{Proof}{#1}\ignorespaces}
    \newcommand{\qed}{\end{Proof}}
    \newcommand{\noqed}{
        \renewcommand{\qedsymbol}{}
        \end{Proof}}}
    \ifthenelse{\boolean{italian}}{
        \renewcommand{\proofname}{Dimostrazione.}}{}

\usepackage[varg]{txfonts}

\usepackage[hypertex]{hyperref}

\begin{document}

\bibliographystyle{amsalpha}

\title{A Factorization Theorem for $\vp$--Maps\thanks{AMS 2000 subject classification 46L08, 46L53, 60G25, 46L55}}
\author{}
\author{
~\\
Michael Skeide\thanks{This work is supported by research funds of University of Molise and Italian MIUR.}
%\sc{Under Construction}\\
}
\date{\vspace{-3.8ex}}
%\date{July 2006}
%$\mtt{\ul{~}}$

{
%\addtolength{\parskip}{-1ex}
\renewcommand{\baselinestretch}{1}
\maketitle

%\vfill
%\newpage

%\vspace{10ex}

%\lf\lf\lf\lf\lf\lf
\begin{abstract}
\noindent
We present a far reaching generalization of a factorization theorem by Bhat, Ramesh, and Sumesh (stated first by Asadi) and furnish a very quick proof.
\end{abstract}

%\tableofcontents
}
%\clearpage

%\vfill

%{\parskip0.5ex plus 0.5ex minus 0.5ex

\section{The result}

Let $E$ and $F$ be Hilbert modules over \nbd{C^*}algebras $\cB$ and $\cC$, respectively. Let $\vp$ be a map from $\cB$ to $\cC$. We say a linear map $T\colon E\rightarrow F$ is a \hl{\nbd{\vp}map} if
\beqn{
\AB{T(x),T(x')}
~=~
\vp(\AB{x,x'})
}\eeqn
for all $x,x'\in E$.

\bthmn
Let $E$ and $F$ be Hilbert modules over unital \nbd{C^*}algebras $\cB$ and $\cC$, respectively. Then for every linear map $T\colon E\rightarrow F$ the following conditions are equivalent:
\begin{enumerate}
\item
$T$ is a \nbd{\vp}map for some  completely positive map $\vp\colon\cB\rightarrow\cC$.

\item
There exists a pair $(\sF,\zeta)$ of a \nbd{C^*}correspondence $\sF$ from $\cB$ to $\cC$ and a vector $\zeta\in\sF$, and there exists an isometry $v\colon E\odot\sF\rightarrow F$ such that
\beqn{
T
~=~
v(\id_E\odot\zeta)
\colon
x
~\longmapsto~
v(x\odot\zeta).
}\eeqn
\end{enumerate}
\ethmn

\proof
$2\Rightarrow 1$.
$\vp:=\AB{\zeta,\bullet\zeta}$ is such a map.

$1\Rightarrow 2$.
By Paschke's \hl{GNS-construction} for CP-maps \cite[Theorem 5.2]{Pas73}, there exist a \nbd{\cB}\nbd{\cC}correspondence $\sF$ and a vector $\zeta\in\sF$ such that $\AB{\zeta,\bullet\zeta}=\vp$ and $\sF=\cls\cB\zeta\cC$. By
\beqn{
\AB{x\odot(b\zeta c),x'\odot(b'\zeta c')}
~=~
c^*\vp(\AB{xb,x'b'})c'
~=~
\AB{T(xb)c,T(x'b')c'},
}\eeqn
$x\odot(b\zeta c)\mapsto T(xb)c$ defines an isometry $v\colon E\odot\sF\rightarrow F$. Specializing to $b=\U_\cB$ and $c=\U_\cC$ we get $v(x\odot\zeta)=T(x)$.\qed

\lf
If $T$ is a \nbd{\vp}map for some completely positive map $\vp\colon\cB\rightarrow\cC$, then the objects in the second part are unique in the following sense.

\bcorn
Suppose $\wt{F}$ is another Hilbert \nbd{\cC}module with a map $\tilde{w}\colon E\rightarrow\wt{F}$ such that the set $\tilde{w}(E)\cC$ is total in $\wt{F}$, and suppose there is an isometry $\tilde{v}\colon\wt{F}\rightarrow F$ such that $\tilde{v}\tilde{w}(x)=T(x)$. Then
\beqn{
u\colon
\tilde{w}(x)
~\longmapsto~
x\odot\zeta
}\eeqn
defines a unitary $\wt{F}\rightarrow E\odot\sF$ where $(\sF,\zeta)$ denotes the (unique) GNS-construction for $\vp$. Moreover, $\tilde{w}$ is a \nbd{\vp}map itself. Alternatively, we may require $\tilde{w}$ to be a \nbd{\vp}map. In that case, $\tilde{v}$ is an isometry, automatically.
\ecorn

\section{Discussion}

The basic ingredient, Paschke's GNS-construction, is a standard result in Hilbert module theory. Paschke's paper \cite{Pas73} and Rieffel's \cite{Rie74} were the first discussing Hilbert modules over not necessarily commutative \nbd{C^*}algebras. So, the GNS-construction is as old as the theory itself. It is a very simple consequence directly from the axioms of Hilbert module and correspondence. (In fact, starting from the definition of completely positive map, the GNS-construction can be used nicely to motivate these axioms.) Another standard ingredient we used, is the tensor product of correspondences. Recall that a \hl{correspondence} from $\cA$ to $\cB$ is a Hilbert \nbd{\cB}module with nondegenerate(\hl{!}) left \nbd{*}action by $\cA$. Recall, too, that the (internal) \hl{tensor product} of a correspondence $\sE$ from $\cA$ to $\cB$ and a correspondence $\sF$ from $\cB$ to $\cC$ is the unique (up to canonical isomorphism) correspondence $\sE\odot\sF$ from $\cA$ to $\cC$ that is generated by elements $x\odot y$ fulfilling
\baln{
\AB{x\odot y,x'\odot y'}
&
~=~
\AB{y,\AB{x,x'}y'},
&
a(x\odot y)
&
~=~
(ax)\odot y.
}\ealn
A Hilbert module $E$ can be considered as a correspondence under the canonical left actions of the algebra of adjointable operators $\sB^a(E)$, the algebra  of compact operators $\sK(E)$ or the complex numbers $\C$.

Another standard result, is inducing a representation of a Hilbert module $E$ by operators in $\sB(G,H)$ from a nondegenerate representation $\pi$ of $\cB$ on a Hilbert space $G$. It helps to recover first the Stinespring construction from GNS-construction and then the result from Bhat, Ramesh, and Sumesh \cite{BRS12} (see the corollary on the next page and the remark following it) from our theorem. This inducing procedure is known since Rieffel's paper \cite{Rie74a} (see the proof of \cite[Proposition 6.10]{Rie74a} in front of the proposition). Indeed, the representation $\pi$ turns $G$ into a correspondence from $\cB$ to $\C$. We define the Hilbert space $H:=E\odot G$. Then every element $x\in E$ gives rise to an operator $L_x\colon g\mapsto x\odot g$ with adjoint $L_x^*\colon y\odot g\mapsto\pi(\AB{x,y})g$. It follows that the map $\eta\colon x\mapsto L_x$  is a representation of $E$ by operators in $\sB(G,H)$ in the sense that $\eta(x)^*\eta(y)=\pi(\AB{x,y})$ and $\eta(xb)=\eta(x)\pi(b)$. Moreover, $H$ is a correspondence from $\sB^a(E)$ to $\C$. In fact, by $\rho(a):=a\odot\id_G$ we define a representation, satisfying $\eta(ax)=\rho(a)\eta(x)$. (Note that also $\eta(x)$ may be written conveniently as $L_x=x\odot\id_G$. The whole induction procedure is nothing but amplifying things with the identity on $G$.)

Note that, given $\pi$, a map $\eta\colon E\rightarrow\sB(G,H')$ (where $H'$ can be any Hilbert space) is uniquely determined (up to suitable unitary equivalence) by the properties $\eta(x)^*\eta(y)=\pi(\AB{x,y})$ and $\eta(xb)=\eta(x)\pi(b)$ plus the cyclicity condition that $\eta(E)G$ be total in $H'$. Note, too, that we need not require that $\eta$ be linear or bounded. (By the stated properties, $\eta$ is a \it{ternary homomorphism} into the Hilbert \nbd{\sB(G)}module $\sB(G,H')$; see Skeide and Abbaspour \cite{AbSk07} for details. In particular, $\eta$ is completely contractive.) Of course, also $\rho$ is determined uniquely by $\eta(ax)=\rho(a)\eta(x)$.

Now let us return to our CP-map $\vp\colon\cB\rightarrow\cC$ with GNS-construction $(\sF,\zeta)$. Suppose we have a representation $\sigma$ of $\cC$ on a Hilbert space $K$. This gives rise to a Hilbert space $G:=\sF\odot K$ and the induced representation $\pi\colon\cB\rightarrow\sB^a(\sF)\rightarrow\sB(G)$ of $\cB$ on $G$. If we put $Z:=L_\zeta\colon k\mapsto\zeta\odot k$, then $Z^*\pi(b)Z=\sigma\circ\vp(b)$. In particular, if $\cC\subset\sB(K)$ is an operator algebra and $\sigma$ is the canonical injection, then $Z^*\pi(b)Z=\vp(b)$. Clearly, $\pi(\cB)ZK=(\cB\zeta)\odot G$ is total in $G$. In other words, $\pi$ is the Stinespring representation \cite{Sti55} of $\cB$ and $Z\in\sB(K,G)$ the cyclic map.

\bcorn
Let $\vp\colon\cB\rightarrow\cC\subset\sB(K)$ be a CP-map and construct its (unique!) Stinespring triple $(G,\pi,Z)$ as explained. Suppose $T\colon E\rightarrow F$ is a \nbd{\vp}map as in the theorem, and let $\eta$ be the (unique!) representation of $E$ into $\sB(G,H)$ induced by the Stinespring representation $\pi$ as discussed. Put $L:=F\odot K$, and define its (unique!) representation $\chi$ into $\sB(K,L)$ induced by the canonical injection $\cC\rightarrow\sB(K)$.

Then with $v$ as in then theorem, $V:=v\odot\id_K\in\sB(E\odot\sF\odot K,F\odot K)=\sB(E\odot G,L)=\sB(H,L)$ is an isometry such that $V\eta(x)Z=\chi\circ T(x)$. Moreover:
\begin{enumerate}
\item
$(G,\pi,Z)$ is determined uniquely by the properties a minimal Stinespring construction fulfills (properties that have nothing to do with the \nbd{\vp}map $T$).

\item
$(H,\eta)$ is determined uniquely by the properties a representation induced by $\pi$ fulfills (properties that have nothing to do with the \nbd{\vp}map $T$).

\item
$V$ is determined uniquely by $V\eta(x)Z=\chi\circ T(x)$. In particular, it is an isometry, automatically.
\end{enumerate} 
\ecorn

\bremn
Like we assumed that $\cC$ is given as a subalgebra of $\sB(K)$ from the beginning, we also may assume that $F$ is a concrete Hilbert \nbd{\cC}module contained in $\sB(K,L)$ from the beginning. (By this we mean that $F$ is a norm closed subspace of $\sB(K,L)$ fulfilling $F\cC\subset F$, $F^*F\subset\cC$, and $\cls FK=L$ .) In this case, $\chi$ reduces to the canonical injection $F\rightarrow\sB(K,L)$ and disappears from the formulae: $V\eta(x)Z=T(x)$. In particular, if $\cC=\sB(K)$ and $F=\sB(K,L)$ we get back the results of \cite{BRS12}.

But a CP-map $\vp\colon\cB\rightarrow\sB(K)$ may be considered as a CP-map into $C^*(\vp(\cB))$ (or any \nbd{C^*}algebra $\cC$ in between), and a \nbd{\vp}map $T\colon E\rightarrow\sB(K,L)$ may be considered a \nbd{\vp}map into the \it{ternary space} $\cls T(E)C^*(\vp(\cB))$ generated by $T(E)$ (or any Hilbert \nbd{\cC}module in between).
\eremn

\bremn
Like the Stinespring construction, also the GNS-construction requires that $\cB$ is unital. The GNS-correspondence can be constructed also when $\cB$ is nonunital, coming shipped with an embedding $i\colon\cB\otimes\cC\rightarrow\sF$ such that $\AB{i(a\otimes b),i(a'\otimes b')}=b^*\vp(a^*a')b'$ and $\sF=\ol{i(\cB\otimes\cC)}$. In order to have a cyclic vector $\zeta$ (without making $\sF$ bigger, what is always possible via unitalization), it is necessary to require that $\vp$ be \it{strict}. For us, the simplest way to describe this, is the condition that for some bounded approximate unit $\bfam{u_\lambda}$ for $\cB$ the corresponding net in $\sB^a(\sF)$ be \nbd{*}strongly convergent to the  identity; see Lance \cite[Sections 2 and 5]{Lan95} and Skeide \cite[Section 4.1]{Ske01} for details.
\eremn

\bremn
The purpose of this note was to illustrate that assuming only little and very basic knowledge about Hilbert modules, results like the Stinespring construction or its generalization to \nbd{\vp}maps drop out easily. (For whom who knows these facts, the corollary in this  section could be stated  already in the end of Section 1. And also there it would not require any proof.) We wish to underline, that this aspect of simplification, though already quite positive as such, is not at all the most important one. Using GNS-construction instead of Stinespring construction has the most striking consequences, when it comes to composition of CP-maps. This is so, because correspondences may be viewed as functors in various ways, and composition of CP-maps, roughly, corresponds to tensor products of correspondences.  Nothing like this works with Stinespring representations! We explained this carefully in Bhat and Skeide \cite[Section 2]{BhSk00}. (In this section the reader can also find a compact introduction to Hilbert modules and correspondences.) We explained it once more very detailedly in the survey \cite{Ske11}. In \cite{BhSk00} the consequent application of GNS-construction instead of Stinespring construction for CP-semigroups led to the first construction of a product system of correspondences.
\eremn

\lf\noindent
\bf{Acknowledgment.~}
I would like to thank the referee for encouraging this revision, putting more effort into the description of the original result of \cite{BRS12}, and for other valuable suggestions.

\newpage

%{\addtolength{\parskip}{-1ex}%\sloppy
\setlength{\baselineskip}{2.5ex}

%\renewcommand{\tt}[1]{\texttt{\small #1}}
% \bibliography{mybib}

\begin{thebibliography}{Rie74b}

\bibitem[AS07]{AbSk07}
G.~Abbaspour and M.~Skeide, \emph{{~~~Generators of dynamical systems on
  Hilbert modules}}, Commun.\ Stoch.\ Anal. \textbf{1} (2007), 193--207,
  (ar\-Xiv: math.OA/0611097).

\bibitem[BRS12]{BRS12}
B.V.R. Bhat, G.~Ramesh, and K.~Sumesh, \emph{{Stinespring's theorem for maps on
  Hilbert $C^*$--modules}}, J.\ Operator Theory \textbf{68} (2012), 173--178,
  (ar\-Xiv: 1001.3743v1).

\bibitem[BS00]{BhSk00}
B.V.R. Bhat and M.~Skeide, \emph{{Tensor product systems of Hilbert modules and
  dilations of completely positive semigroups}}, Infin.\ Dimens.\ Anal.\
  Quantum Probab.\ Relat.\ Top. \textbf{3} (2000), 519--575, (Rome,
  Volterra-Pre\-print 1999/0370).

\bibitem[Lan95]{Lan95}
E.C. Lance, \emph{{Hilbert $C^*$--modules}}, Cambridge University Press, 1995.

\bibitem[Pas73]{Pas73}
W.L. Paschke, \emph{{Inner product modules over $B^*$--algebras}}, Trans.\
  Amer.\ Math.\ Soc. \textbf{182} (1973), 443--468.

\bibitem[Rie74a]{Rie74}
M.A. Rieffel, \emph{{Induced representations of $C^*$-algebras}}, Adv.\ Math.
  \textbf{13} (1974), 176--257.

\bibitem[Rie74b]{Rie74a}
\bysame, \emph{{Morita equivalence for $C^*$--algebras and $W^*$--algebras}},
  J.\ Pure Appl.\ Algebra \textbf{5} (1974), 51--96.

\bibitem[Ske01]{Ske01}
M.~Skeide, \emph{{Hilbert modules and applications in quantum probability}},
  Ha\-bi\-li\-ta\-tions\-schrift, Cottbus, 2001, Available at
  {\footnotesize\url{ http://web.unimol.it/skeide/}}.

\bibitem[Ske11]{Ske11}
\bysame, \emph{{Hilbert modules---square roots of positive maps}}, Quantum
  Probability and Releated Topics --- Proceedings of the 30th Conference
  (R.~Rebolledo and M.~Orszag, eds.), Quantum Probability and White Noise
  Analysis, no. XXVII, World Scientific, 2011, (ar\-Xiv: 1007.0113v1),
  pp.~296--322.

\bibitem[Sti55]{Sti55}
W.F. Stinespring, \emph{{Positive functions on $C^*$--algebras}}, Proc.\ Amer.\
  Math.\ Soc. \textbf{6} (1955), 211--216.

\end{thebibliography}

\newcommand{\Swap}[2]{#2#1}\newcommand{\Sort}[1]{}
\providecommand{\bysame}{\leavevmode\hbox to3em{\hrulefill}\thinspace}
\providecommand{\MR}{\relax\ifhmode\unskip\space\fi MR }
\providecommand{\MRhref}[2]{%
  \href{http://www.ams.org/mathscinet-getitem?mr=#1}{#2}
}
\providecommand{\href}[2]{#2}

\noindent
Michael Skeide,
{\small\itshape Universit\`a\ degli Studi del Molise},
{\small\itshape Dipartimento S.E.G.e S.},
{\small\itshape Via de Sanctis}
{\small\itshape 86100 Campobasso, Italy},
{\small{\itshape E-mail: \tt{skeide@unimol.it}}},\\
{\small{\itshape Homepage: \tt{http://www.math.tu-cottbus.de/INSTITUT/lswas/\_skeide.html}}}.

%\listofOWs

\end{document}